\documentclass[a4paper,12pt]{article}

\pdfoutput=1

\usepackage[a4paper,includeheadfoot,top=25mm,bottom=25mm,footskip=7mm,pdftex]{geometry}  

\usepackage[T1]{fontenc}
\usepackage[latin1]{inputenc}  %% Buchstaben ä,ö,ü,ß
\usepackage[english]{babel}    %% Silbentrennung und erkennen von "a,"o,"u,"s
\usepackage{amsmath,amssymb,amsfonts,textcomp}
\usepackage{graphicx}
\usepackage{xcolor}
\newcommand{\C}{\mathbb{C}}			% die komplexen Zahlen
\newcommand{\R}{\mathbb{R}}			% die reelen Zahlen
\newcommand{\N}{\mathbb{N}}			% die natürlichen Zahlen 
\renewcommand{\P}[1]{\mathbb{P}(#1)}			% Partitionen
\newcommand{\eins}{\boldsymbol{1}}			% Einselement in Algebren
\newcommand{\expo}[2]{\mathrm{e}_{#1 }{}^{#2 }} %Faltungsexponential
\newcommand{\bigop}[3]{\overset{#3}{\underset{#2}{\pmb{#1}\,}}} %
\newcommand{\pa}[1]{\left( #1 \right)} % variable runde Klammern
\newcommand{\pb}[1]{\left[ #1 \right]} % variable eckige Klammern
\newcommand{\pd}[1]{\left< #1 \right>} % variable spitze Klammern
\newcommand{\norm}[1]{\left\| #1 \right\|} % variable doppel striche Klammern
\newcommand{\nbd}{\nobreakdash-\hspace{0pt}} % Kein Zeilenumbruch an Bindestrich
\newcommand{\astalg}{\star} % Convolution of Alg-Hom's over dual semi-groups

\usepackage{hyperref}
\newtheorem{lm}{Lemma}[section]

\newtheorem{sz}[lm]{Theorem}
\newtheorem{fo}[lm]{Corollary} 
\newtheorem{df}[lm]{Definition} 
\newtheorem{bm0}[lm]{Remark}

\newtheorem{bsp}[lm]{Example}

\newenvironment{beweis}{{\em Proof:}}{\hfill $\square$\goodbreak \bigskip}

\begin{document}
\author{Stefan Vo\ss \\
Department of Mathematics and Computer Science  \\
University of Greifswald  \\
Walther-Rathenau-Stra\ss e 47  \\
17487 Greifswald, Germany   \\
}
\title{Unitary cocycles and processes on \\the full Fock space}
\maketitle
\begin{abstract}
We consider a unitary cocycle or Sch\"urmann triple on the non-com\-mu\-ta\-tive unitary group fixed by a complex matrix which induces an additive free white noise or an additive free L\'evy process on the tensor algebra over the full Fock space. A L\'evy process on a Voiculescu dual semi-group is given by a generator or Sch\"urmann triple. We will show how a free L\'evy process on the non-commutative unitary group fixed by a complex matrix  can be obtained by infinitesimally convolving the additive free white noise. 
%We present a theory of how every free L\'evy process on the non-commutative unitary group fixed by a complex matrix can be described by additive free white noise by the method of infinitesimal convolution.  
\end{abstract}

\pagenumbering{arabic}

\newcommand{\KD}{\operatorname{K}\langle d \rangle} % 

\section{Introduction}
From an algebraic point of view, a (stochastic) process is a family $(f_{t})_{t}$ of unital $*$\nbd algebra homomorphisms $f_{t}: B \to A$ on a so called quantum probability space (QPS). This is a pair $(A,\Phi)$ consisting of a unital $*$\nbd algebra $A$ and a state $\Phi$, that is a positive linear normed functional $\Phi:A\to \C$ and playing the role of an expectation. A well known example is the space of all linear adjointable maps of the symmetric or Boson Fock space \cite{Part92} where $\Phi$ denotes the vacuum expectation. 
Another example is the space of all linear adjointable maps of the full Fock space $\Gamma(H)$ over the Hilbert space $H$ and the state $\Phi(g):=\langle \Omega, g(\Omega)\rangle$ for a linear adjointable map $g:\Gamma(H) \to \Gamma(H)$. A thorough introduction to quantum probability can be found in \cite{Part92,Hol01}.

From a measure-theoretical point of view, a stochastic process $X$ is a family of measurable maps $X_{t}:E\to G$, where $E$ is a probability space. Given a stochastic process $X=(X_{t})_{t\in \R}$ we get an algebraic process $f_{t}:B\to A$ by $f_{t}(\varphi):=\varphi  \circ X_{t}$, where $B:=L^{\infty}(G)$ and $A:=L^{\infty}(E)$. Since $L^{\infty}(G)$ is a commutative unital algebra we would like to stress that $B$ may well be a non-commutative algebra, sometimes it carries additional structure making it a quantum group, a Hopf algebra or in this paper a dual group in the sense of Voiculescu \cite{Voi1,Voi87}. For example, the authors of \cite{Lindsay12} investigate quantum groups and their analytic aspects in a series of papers. Essentially, the additional structure for dual groups in this paper will be an associative convolution $\astalg$ of unital $*$\nbd algebra homomorphisms and a unit element $\delta$ with respect to this convolution.

The algebraic version of factorisation also known as independence, i.e.\ the joined distribution equals the product of the marginal distributions,   for unital $*$\nbd algebra homomorphisms $f$ and $g$ is $\Phi \circ (f \sqcup g) = (\Phi \circ f)\odot (\Phi \circ g)$ . In contrast to the independence in (classical) probability theory there is more than one way to choose the product $\odot$, see the 'Muraki five': tensor, free, boolean, monotone and anti-monotone product \cite{MurN03}. 
% \cite{MurN01,MurN02,MurN03}. 

A L\'evy process should have independent and stationary increments. So from our algebraic point of view, see also \cite{MSchue93,BGhSc,Fra06}, a L\'evy process is a family $(f_{s,t})_{s\leq t}$ of unital $*$\nbd algebra homomorphisms $f_{s,t}:B\to (A,\Phi)$ which are called increments, i.e.\ factorise with respect to the convolution $f_{r,t} = f_{r,s}\astalg f_{s,t}$ for all $r\leq s \leq t$. The increments of different time periods factorise with respect to $\odot$. The expectation of an increment $f_{s,t}$ is stationary, i.e.\ only depends on $(t-s)$. Moreover, it is (weakly) continuous, i.e.\ $\Phi \circ f_{s,t}(b)$ converges to $\delta(b)$ when $t$ tends to $s$. We investigate  free L\'evy processes on the non-commutative unitary group $\KD $ over the full Fock space. This investigation also yields a family $(U_{r,s})_{r\leq s}$ of unitaries, which is indeed a process with independent increments. In this paper, $\KD $ is only considered for finite $d$. What happens when $d$ is infinite can be found in \cite{Sinha10a}. The non-commutative unitary group $\KD $ also introduced in \cite{vW84} is the algebra of all polynomials of $d\times d$ non-commuting indeterminates $x_{k,l},x_{k,l}^{*}$ with unitary relations $X^{*}X-E$, $XX^{*}-E$, where $X$ is the matrix with entries $x_{k,l}$ and $E$ is the $d\times d$ identity matrix. The convolution is defined by $f\astalg g (x_{k,l}):=f \sqcup g ((X^{(1)} \cdot X^{(2)})_{k,l})$ and the counit is the Kronecker delta  $\delta(x_{k,l})=\delta_{k,l}$.

The notion of a cocycle appears on various occasions \cite{Part92}%\cite{McL75}
. In this paper, we give Sch\"urmann's algebraic definition of a unitary cocycle, also called Sch\"urman triple $(\rho,\eta,\Psi)$. We state the proper definition in Section \ref{def:uc} where (ii) is a reformulation of a cocycle condition. 

We will infinitesimally convolve the additive free white noise $T(I_{s,t})$, see Theorem \ref{addWhiteNoise}, given by a Sch\"urmann triple $(\rho,\eta,\Psi)$ coming from a complex $d\times d$ matrix $L$ as in Theorem \ref{kdGenerator} with respect to the convolution of the non-commutative unitary group $\KD $, show that it is a convergent net and that the limit forms a $\odot$-free L\'evy process on $\KD $ over the full Fock space with generator $\Psi$. Finally, infinitesimally convolving the limit by the additive convolution ensures that this L\'evy process is cyclic, that is applying the L\'evy process repeatedly on the vacuum vector $\Omega$ yields a dense subspace of the full Fock space. It follows that our L\'evy process with generator $\Psi$ over the full Fock space is stochastically equivalent to every $\odot$-free L\'evy process with the same generator $\Psi$ on $\KD $ over any QPS.

A general framework accommodating any dual semi-group, the 'Muraki five' products and any Sch\"urmann triple where a representation theorem over the respective Fock space is aquired by applying the same method is in preparation by the author.

% A general framework of this method of infinitesimal convolution for every dual semi-group, the 'Muraki five' and for every Sch\"urmann triple and so yielding the representation theorem on Fock spaces is in preparation by the author. 

On the level of Hopf algebras, this method of infinitesimal convolution is introduced in \cite{ScSV10}.

In Section \ref{sec:Pre}, we state the necessary definitions and prove  convergence to a 'convolution exponential' for dual semigroups. In Section \ref{sec:lem}, we present a main lemma for normed unital algebras which is used to infinitesimally convolve the free additive white noise in Section \ref{sec:example} on the non-commutative unitary group in \ref{sec:ugroup} and on the tensor algebra in \ref{sec:tensoralg}. In the final section we construct a free L\'evy processes by infinitesimal convolution.

\section{Preliminaries}\label{sec:Pre}
We will consider associative algebras over $\C$, the field of complex numbers.
A $*$-algebra is an algebra $A$ with an involution $*$, i.e.\ an anti-linear map
$a \mapsto a ^{*}$ on $A$ such that $(ab) ^{*} = b ^* a ^*$ and $(a ^* )^* = a$.
A unital algebra is an algebra containing an element $\eins$ (called unit element) in $A$ with
$a \, \eins = a = \eins \, a$.\par
A partition $\alpha$ is a finite subset of a closed interval $[ S,T ]\subset \R$. It consists of elements $S=t_{1}<\dotsb < t_{n}=T$ and is denoted by $\left\{  S=t_{1}<\dotsb < t_{n}=T \right\}$. Let $\P{[ S,T ]}$ be the set of all partitions of $[ S,T ]$ partially ordered by inclusion. For a family of maps $(k_{t})_{t}$ indexed by one parameter $t$ we write 
\begin{align*}
\sum\limits_{\alpha} k_{\alpha}:= k_{t_{2}-t_{1}}+\dotsb + k_{t_{n}-t_{n-1}}.
\end{align*}
For a family of maps $(k_{s,t})_{s,t}$ indexed by two parameters $s$ and $,t$ we write 
\begin{align*}
\sum\limits_{\alpha} k_{\alpha}:= k_{t_{1},t_{2}}+\dotsb + k_{t_{n-1},t_{n}}. 
\end{align*}
We use the same notation for other binary operators as well. A family of elements $(\theta_{\alpha})_{\alpha \in \P{[ S,T ]}}$ is a convergence net to $\theta$ if $\theta = \lim\limits_{\alpha\in \P{[ S,T ]}}\theta_{\alpha} $, i.e. for all $\epsilon>0$ there exists a partition $\gamma\in \P{[S,T]}$ such that $\norm{ \theta_{\alpha} -\theta} <\epsilon$ for all $\alpha \geq \gamma$. The family $(\theta_{\alpha})_{\alpha \in \P{[ S,T ]}}$ is a Cauchy net if for all $\epsilon>0$ there exists a partition $\gamma\in \P{[S,T]}$ such that $\norm{ \theta_{\alpha} -\theta_{\beta} } <\epsilon$ for all $\alpha \geq \gamma ,\beta \geq \gamma$. Similar to the theory of sequences, a metric space is complete if and only if all Cauchy nets are convergence nets.\par
We write $A \sqcup B$ for the free product of the algebras $A$ and $B$, and $A \sqcup_{1} B$ for the free product of unital algebras. Moreover, let $\iota_{A}$ (resp.\ $\iota_{B}$) be the canonical embedding. We often identify $A$ (resp.\ $B$) with its embedding in $A \sqcup B$ or $A \sqcup_{1} B$. Either one is the co-product in its respective category, namely the categories of algebras and unital algebras \cite[Section 2]{BGhSc}. Let $k_{1} \amalg k_{2}: A_{1} \sqcup A_{2} \to B_{1} \sqcup B_{2}$ be the algebra homomorphism defined by $k_{1} \amalg k_{2}:= (\iota_{1} \circ k_{1}) \sqcup (\iota_{2} \circ k_{2})$ and  $k_{1} \amalg_{1} k_{2}$ analogously. In the category of pairs formed by  $(B,\varphi)$, where $B$ is an algebra and $\varphi:B\to \C$ a linear functional  we will use an abstract notion of independence. 
\begin{df}\label{PUP} The notion of independence is given by a tensor structure of the form
\begin{align*}
(B_{1},\varphi_{1})\Box (B_{2},\varphi_{2}) = (B_{1}  \sqcup B_{2}, \varphi_{1}\odot \varphi_{2})
\end{align*}
in the tensor category with inclusions, cf. \cite{Fra06}. This means that $\odot$ satisfies the following conditions: 
The map $\varphi_{1} \odot \varphi_{2}$ is a linear functional and 
\begin{align*}
(\varphi_{1} \odot \varphi_{2}) \circ \iota_{i} &= \varphi_{i} \quad i=1,2\tag{UP1}\\
(\varphi_{1} \odot \varphi_{2})\odot \varphi_{3} &= \varphi_{1}\odot (\varphi_{2} \odot \varphi_{3})\tag{UP2}\\
(\varphi_{1} \circ k_{1})\odot (\varphi_{2} \circ k_{2}) &= (\varphi_{1} \odot \varphi_{2}) \circ (k_{1}  \amalg k_{2})\tag{UP3}
\end{align*}
for linear functionals $\varphi_{i}:B_{i} \to \C$  $(i=1,2,3)$ algebra homomorphism $k_{i}:A_{i} \to B_{i}$, $(i=1,2)$ and for all $b_{1}\in B_{1}$ , $b_{2}\in B_{2}$.
\end{df}
If  we assert
\begin{align*}
(\varphi_{1} \odot \varphi_{2})(b_{1}b_{2}) &= \varphi_{1} \odot \varphi_{2}(b_{2}b_{1})=\varphi_{1}(b_{1}) \cdot \varphi_{2}(b_{2})\tag{UP4}
\end{align*}
or for $*$-algebras equivalently 
\begin{align*}
\widetilde{\varphi_{1}}, \widetilde{\varphi_{2}} \text{ states } \Rightarrow \widetilde{\varphi_{1}\odot \varphi_{2}} \text{ state.}\tag{UP4'}
\end{align*}
then there are exactly five examples of $\odot$, namely the Muraki five: tensor, boolean, free, (anti-)monotone \cite{MurN03}. % \cite{MurN01,MurN02,MurN03}. 
 A brief overview is available in \cite[Section 2]{BGhSc}. We will only discuss the free product also known as freeness \cite{NS06}. There are two ways to describe the free product: Let  $b_{1}, \ldots, b_{m}$ be in $B_{1}$ or $B_{2}$  ($\subset B_{1} \sqcup B_{2}$) where consecutive elements are contained in different algebras
\begin{align}
\varphi_{1}\odot \varphi_{2}\left( b_{1}\dotsb b_{m} \right) := \sum\limits_{\makebox[0pt][c]{$\mathsurround=0pt\scriptstyle{ I\underset{\not=}{\subset}\left( 1,\ldots,m \right)}$} }( -1 )^{m-\#  I+1}\varphi_{1}\odot \varphi_{2}\biggl( \prod\limits_{k\in I}^{\to}b_{k} \biggr)  \prod\limits_{j\notin I}\varphi_{1}\odot \varphi_{2}\left( b_{j} \right) \label{fp::A}
\end{align}
% \begin{align}
% \varphi_{1}\odot \varphi_{2}\left( b_{1}\dotsb b_{m} \right)\! :=\!\! \sum\limits_{I\underset{\not=}{\subset}\left( 1,\ldots,m \right) }\left( -1 \right)^{m-\#  I+1}\!\varphi_{1}\!\odot\! \varphi_{2}\!\biggl( \prod\limits_{k\in I}^{\to}b_{k} \biggr)\!  \prod\limits_{j\notin I}\varphi_{1}\odot \varphi_{2}\left( b_{j} \right) \label{fp::A}
% \end{align}
as recursion formula with $\varphi_{1}\odot \varphi_{2}\left( \prod_{k\in \emptyset}^{\to}b_{k} \right):=1 \in\C$. For normed linear maps $\varphi_{1},\varphi_{2}$ of unital algebras $B_{1},B_{2}$ we have
\begin{align}
\varphi_{1}\odot \varphi_{2} (b_{1}\dotsb b_{m})=0 \quad \text{if } (\varphi_{1}\odot \varphi_{2}) (b_{k})=0\quad \text{for all } k=1, \ldots, m. \label{fp::B}
\end{align}
For unital algebras together with normed linear maps both (\ref{fp::A}) and (\ref{fp::B}) are equivalent.
Elaborating on (\ref{fp::A}), consider the example $\prod_{k\in \left\{ 1,5,6,12, 19  \right\}}^{\to}b_{k}=b_{1}b_{5}b_{6}b_{12}b_{19}$ where the product is in $B_{1} \sqcup B_{2}$ and assume $b_{1}\in B_{1}$. It is important to be aware that $b_{1}b_{5}$ is one element in the next step of the recursion because $b_{1}b_{5}$ is a product within the algebra $B_{1}$.\par
Two algebra homomorphisms $f:B_{1} \to (A,\Phi)$ and $g:B_{2} \to (A,\Phi)$ are $\odot$-independent if $\Phi \circ(f \sqcup g) = (\Phi  \circ f) \odot (\Phi \circ g)$.

We will reformulate the definition of a dual group given in the original paper by Voiculescu \cite{Voi1} in the following way: 
\begin{df}A dual semigroup $\left( B, \Delta,\delta \right)$ consists of a unital $*$\nbd algebra $B$ and unital $*$\nbd algebra homomorphisms $\Delta:B\to B  \sqcup_{1} B$, $\delta:B\to \C$ such that
\begin{gather}
\left( \Delta\amalg_{1} id_{B} \right) \circ \Delta = \left( id_{B}\amalg_{1} \Delta \right) \circ \Delta \\
\left( \delta \amalg_{1} id_{B} \right) \circ \Delta = id_{B} 
= \left( id_{B} \amalg_{1} \delta \right) \circ \Delta 
\end{gather}
holds. If, in addition, there exists a unital $*$\nbd algebra homomorphism $S':B\to B$ such that $\left( S'  \sqcup_{1} id_{B} \right) \circ \Delta = \delta
= \left( id_{B}  \sqcup_{1} S' \right) \circ \Delta$ then $B$ is called a dual group.
%, where $ \delta \amalg_{1} id_{B} $ und $id_{B} \amalg_{1} \delta $ is to be  understood as projections of $B \sqcup_{1} B$ onto the left resp. right $B$. 
\end{df}
\begin{bsp}[Tensor algebra]\label{dSG:tensoralgebra} Let $V$ be a $*$\nbd vector space, i.e.\ there exists an anti-linear map $v\mapsto v^{*}$ with $(v^{*})^{*}=v$. The tensor algebra is the direct sum $T(V):=\C\eins  \oplus V \oplus V^{ \otimes 2} \oplus \dotsb$ with multiplication $v \cdot w := v \otimes w$ possessing the universal property: For a unital $*$\nbd algebra $A$ and a linear map $R:T(V)\to A$ there exists a unique unital $*$\nbd algebra homomorphism $T(R):T(V) \to A$ with $T(R) \restriction V = R$. Furthermore, $T(R)$ is given by $T(R)(\eins)=\eins$ and \goodbreak $T(R)(v_{1} \otimes \dotsb  \otimes v_{n})=(R(v_{1})) \dotsb (R(v_{n}))$%\cite[III§5.1 Prop.1]{BA}%Bourbaki
.\par
The tensor algebra is a dual semigroup with $ \Delta:= T\left( f \right)$ and $\delta:= T\left( 0 \right)$, where $
f :V \to T\left( V \right) \sqcup_{1} T\left( V \right)$ with $f\left( v \right) = v^{\left( 1 \right)}+v^{\left( 2 \right)}$. The vector $v^{\left( 1 \right)}$ (resp.\ $v^{\left( 2 \right)}$) is in the left (resp.\ right) $V$ component of $T\left( V \right) \sqcup_{1} T\left( V \right)$. Together with the extension of the map $v\mapsto -v$ to $S'$ by the universal property the tensor algebra is a dual group.
\end{bsp}
\begin{bsp}[Non-commutative unitary group]\label{dSG:kd}
Let $\KD$ be the unital \emergencystretch 1.0em $*$\nbd algebra of polynomials in non-commuting indeterminates $x_{kl}$ and $x_{kl}^{*}$ for $k,l=1,\ldots,d$ with complex coefficients and relations
% \begin{align*}
% \sum\limits_{n=1}^{d} x_{kn}x_{ln}^{*} - \delta_{kl}\eins{}, \quad\quad \sum\limits_{n=1}^{d} x_{nk}^{*}x_{nl} -\delta_{kl}\eins{},
% \end{align*}
%
\begin{align*}
\left( X \cdot X^{*} - E \right)_{k,l},\quad \quad \left( X^{*} \cdot X -E \right)_{k,l}\quad \forall k,l=1,\ldots,d
\end{align*}
where $X_{k,l}:=x_{k,l}$, $X_{k,l}^{*}:=x_{l,k}^{*}$ and $E$ the $d\times d$ identity matrix
 and involution $x_{kl} \mapsto x_{kl}^{*}$. 
A unital $*$\nbd algebra homomorphism on $\C\langle d \rangle$ (the above without relations)  is determined by the values of  $x_{kl}$ for $k,l=1,\ldots,d$%, see \cite[III§2.6 Prop.6,§2.7 Prop.7]{BA} % Bourbaki
.

Define $\Delta: \C\langle d\rangle \to \KD  \sqcup_{1} \KD $ and $\delta:\C\langle d\rangle \to \C$ by
% \begin{align*}
% \Delta \left( x_{kl} \right)= \sum\limits_{n=1}^{d} x_{kn}^{\left( 1 \right)} x_{nl}^{\left( 2 \right)}, \quad\quad \delta\left( x_{kl} \right) = \delta_{kl}
% \end{align*}
\begin{align*}
\Delta \left( x_{kl} \right)= \left( X^{\left( 1 \right)} X^{\left( 2 \right)} \right)_{k,l}, \quad\quad \delta\left( x_{kl} \right) = E_{kl}
\end{align*}
for all $k,l$. Since the maps $\Delta$ and $\delta$ respect the relations above, $\left(\KD , \Delta,\delta\right)$ is a dual semigroup. Moreover, it is a dual group with $S'(x_{k,l}):=x_{l,k}^{*}$.
\end{bsp}
A convolution of unital $*$\nbd algebra homomorphisms $f$ and $g$ and a convolution of functionals $\varphi_{1}$ and $\varphi_{2}$ is defined by 
\begin{align*}
f\astalg  g:= (f    \sqcup_{1} g)  \circ \Delta && \varphi_{1} \uplus \varphi_{2}:= (\varphi_{1} \odot \varphi_{2}) \circ \Delta.
\end{align*}
The proof of the following theorem transfers the convolution $\uplus$ to the coalgebra convolution in the symmetric tensor algebra $(S(V),\ast)$, where $\ast$ is determined by the choice of $\odot$ and $\Delta$ \cite[Theorem 3.4]{BGhSc}. For a linear map $\Psi:V\to \C$ we define $D(\Psi):S(V)\to\C$ by $D(\Psi) \restriction V = \Psi$ and zero elsewhere. Due to the fundamental theorem of coalgebras \cite{DNR} the exponential 
\begin{align*}
\exp_{\ast}(tD(\Psi)) := \sum\limits_{n=0}^{\infty} \frac{t^{n}D(\Psi)^{\ast n} } {n!}= S(0) + tD(\Psi) +\frac{t^{2}}{2}D(\Psi)^{\ast 2} + \dotsb
\end{align*}
 exists point-wise %\cite{MSchue85},
\cite[Section 4]{MSchue90a}%\cite[page 160]{BGhSc}.
\begin{sz}\label{dualGroup:expo}
Let $(k_{t})_{t\geq 0}$ be a family of functionals given by
\begin{align*}
k_{t}:= \delta + t\Psi + R_{t} \quad \forall t\geq 0
\end{align*}
with $\Psi(\eins)=0$, $R_{t}$ linear, $R_{t}(\eins)=0$ and assume that for every $b\in B$ there exist $C_{b},\epsilon_{b} \geq 0$ with $|R_{t}(b)|\leq t^{2}C_{b}$ for all $t\leq \epsilon_{b}$. Then
\begin{align*}
\lim\limits_{\alpha\in \P{[ S,T ]}} \bigop{ \uplus}{\alpha}{}k_{\alpha}(b) = \quad \exp_{\ast}((T-S)D(\Psi)) (b).
\end{align*}
\end{sz}
\begin{beweis}
Recall that $S(f \uplus g)=S(f)\ast S(g)$ \cite[Theorem 3.4]{BGhSc}. For 
$s\in S(B_{0})$ the following holds
\begin{align*}
 S\left( k_{t} \right)(s)=\left( S(0) +t D (\Psi ) + N_{t}\right)(s)=: \hat{k_{t}}(s)
\end{align*}
for $N_{t}$ like $R_{t}$. This is proven by induction over $s\in B_{0}^{ \otimes_{s} n}$. 
Then
\begin{align*}
 \bigop{ \uplus}{\alpha}{}k_{\alpha}(b)&=S\left( \bigop{ \uplus}{\alpha}{}k_{\alpha} \right)(b) = \bigop{ \ast}{\alpha}{} S\left( k_{\alpha} \right)(b)\\
 &=\left(\bigop{ \ast}{\alpha}{} \hat{k}_{\alpha}\right)(b) \xrightarrow[\alpha\in \P{[ S,T ]}  ]{\text{net}}  \exp_{\ast}((T-S)D(\Psi)) (b)
\end{align*}
convergences by \cite[Lemma 4.2]{ScSV10}. See also \cite[Lemma 3.2]{SV12}
\end{beweis}
The full Fock space $\Gamma(H)$ over a complex Hilbert space $H$ is the direct sum of all tensor powers $H^{ \otimes n}$ of $H$, $n\geq 0$, where $H^{ \otimes 0}:=\C$. When taking the symmetric tensor powers we obtain the well known Fock space or Boson Fock space \cite{Part92}. Define the vacuum vector $\Omega:=(1,0,0, \ldots) \in \Gamma(H)$. 
%$\eins \in \C$ and zero for all other tensor powers of $H$. 
Let $D$ be a pre-Hilbert space and let $L_{a}(D)$ be the unital $*$\nbd algebra of all linear adjointable maps from $D$ to $D$.  We define the creation and annihilation operators $A^{*}, A:D \to L_{a}(\Gamma(H))$ and the preservation operator $\Lambda:L_{a}(D)\to L_{a}(\Gamma(H))$ on a dense subspace $D$ of $H$ by
\begin{itemize}
\item $A^{*}(d)\Omega:=d$ and $A^{*}(d)(d_{1} \otimes \dotsb  \otimes d_{n}):=d \otimes d_{1} \otimes \dotsb  \otimes d_{n}$,
\item $A(d)\Omega:=0 \in \Gamma(H)$ and $A(d)(d_{1} \otimes \dotsb  \otimes d_{n}):=\langle d ,d_{1}\rangle (d_{2} \otimes \dotsb  \otimes d_{n})$,
\item $\Lambda(T)\Omega:=0\in \Gamma(H)$ and $\Lambda(T)(d_{1} \otimes \dotsb  \otimes d_{n}):= (T(d_{1})) \otimes \dotsb  \otimes d_{n}$.
\end{itemize}
The creation and annihilation operator are the adjoints of each other and have norm $\norm{A^{*}(h)}=\norm{A(h)}=\norm{h}$. See also \cite{BoKuSp96}.

\section{Statement of main lemma}\label{sec:lem}
We state the main lemma for unital normed algebras. In order to keep the proof of the main Lemma \ref{abstr:sz:cauchyNet} simple we extract Lemma \ref{BehII}.
%To shorten the proof of the Theorem \ref{abstr:sz:cauchyNet} we extract a starting lemma and then state the theorem and prove it.
%\subsection{Lemma}
%
\begin{lm}\label{BehII}
Let $A$ be a unital normed algebra and $0\leq S< T<\infty$. If there are a constant $C>0$, $ S=s_{1} < \dotsb < s_{n+1}=T$ and $a_{1},\ldots ,a_{n} \in A$ for some  $n\in \N$ such that  $\norm{ a_{i}} \leq (s_{i+1}-s_{i})C$ for all $i=1\ldots n$. Then
\begin{align*}
\norm{\prod\limits_{i=1}^{n} \left( 1+a_{i} \right)\quad -1 - \sum\limits_{i=1}^{n} a_{i} } 
\, \leq \,\,    \frac{1}{2}(T-S)^2C^2 \expo{}{C(T-S)}{}.
\end{align*}

\end{lm}
\begin{beweis}
In order to reshape the product define $c_{1}^{(j)}:= 1,
c_{2}^{(j)}:= a_{j}$ for $j=1 \ldots n$. This yields 
$
\prod\limits_{i=1}^{n} \left( 1+a_{i} \right) = \sum\limits_{k=(k_{1}, \ldots,k_{n}) \in \{1,2\}^{n} } \prod\limits_{j=1}^{n} c_{k_{j}}^{(j)}
$.
Let $D$ be the set  
\begin{align*}
\{1,2\}^{n} \setminus \bigl( \{(1, \ldots , 1)\} 
\cup  \{ (1, \ldots,1,\underset{j\text{-th}}{2},1, \ldots,1), j = 1  \ldots n \} \bigr).
\end{align*}
Then
\begin{align*}
\sum\limits_{k=(k_{1}, \ldots,k_{n}) \in \{1,2\}^{n} } \prod\limits_{j=1}^{n} c_{k_{j}}^{(j)} -1 - \sum\limits_{i=1}^{n} a_{i}  = \sum\limits_{k=(k_{1}, \ldots,k_{n}) \in D } \prod\limits_{j=1}^{n} c_{k_{j}}^{(j)} 
\end{align*}
since the summand equals $1$ for $k=(1,\ldots,1)$  and $a_{j}$ for $k = (\ldots 1,\underset{j\text{-th}}{2},1 \ldots)$.
We will now prove the estimation stated above. Define $Z:= \sum\limits_{i=1}^{n}||a_{i}||$. The constraints on $a_{i}$ imply that $Z\leq (T-S)C$. Therefore,
\begin{align*}
&\norm{\sum\limits_{k=(k_{1}, \ldots,k_{n}) \in D } \prod\limits_{j=1}^{n} c_{k_{j}}^{(j)}} 
\leq \sum\limits_{i=2}^{n} \,\,\sum\limits_{1\leq k_{1} < \dotsb < k_{i}\leq n} \,\,\prod\limits_{j=1}^{i} \norm{a_{k_{j}}}\\
&= (1+\norm{a_{1}})  \dotsb   (1+\norm{a_{n}}) - (\norm{a_{1}} + \dotsb + \norm{a_{n}}) - 1 \\
&\leq \expo{}{Z}{} - Z -1 = \frac{Z^2}{2}(1+ \frac{2Z}{3!}+\frac{2 Z^2}{4!}+\dotsb) \leq \frac {Z^2}{2}\expo{}{Z}{}\\
&\leq  \frac{C^2}{2}(T-S)^2 \expo{}{(T-S)C}{}
\end{align*}
since $(1+||a_{1}||)  \dotsb   (1+||a_{n}||) \leq \expo{}{||a_{1}||}{}\dotsb \expo{}{||a_{n}||}{} = \expo{}{Z}{}$ and $\frac{2}{(n+2)!}\leq \frac{1}{n!}$ for all $n\geq 0$.
\end{beweis}
%
%
%
%
%\subsection{Theorem}
We proceed with the main lemma. Therefore, consider a unital normed algebra $A$ and constants $C>0$, $R<S \in \R_{+}$. Let  $\left( g_{r,s} \right)_{r,s} \in A$  be a family with 
\begin{align*}
g_{r,s}:= 1+a_{r,s} \quad \forall r\leq s \in [ R,S ]
\end{align*}
 satisfying
\begin{align}
\norm{ a_{r,s}}&\leq (s-r)C  \label{PCondC}\\
 r<s<t \quad &\Rightarrow \quad  a_{r,t}= a_{r,s}+a_{s,t}.\label{PCondE}
\end{align}
We observe that
\begin{align*}
&\norm{g_{r,s}}\leq \left( 1+(s-r)C \right) \leq \expo{}{(s-r)C}{}\\
&\norm{g_{t_{1},t_{2}}}\dotsb \norm{g_{t_{n},t_{n+1}}} \leq \expo{}{(t_{n+1}-t_{1})C}{}.
\end{align*}
\begin{df} For $\alpha = \left\{ R=t_{1}<\dotsc < t_{n+1} =S  \right\} \in \P{[R,S]}$
we define
\begin{align*}
\Theta_{\alpha} := (g_{t_{1},t_{2}}) \dotsb (g_{t_{n},t_{n+1}}).
\end{align*}
\end{df}
The norm of $\Theta_{\alpha}$ is bounded by $\expo{}{(S-R)C}{}$.
\begin{lm}\label{abstr:sz:cauchyNet}
The net $(\Theta_{\alpha})_{\alpha \in \P{[R,S]}}$ is a Cauchy net.
\end{lm}
%
%}
\begin{beweis}
The lemma holds if for all $\epsilon>0$ there exists a partition $\gamma\in \P{[R,S]}$ such that $\norm{ \Theta_{\alpha} -\Theta_{\beta} } <\epsilon$ for all $\alpha \geq \gamma ,\beta \geq \gamma$.
Let $\alpha,\beta,\gamma \in \P{[R,S]}$ with $\alpha \geq \gamma ,\beta \geq \gamma$. Then
\begin{align*}
\norm{ \Theta_{\alpha} -\Theta_{\beta} } &= \norm{ \Theta_{\alpha} -\Theta_{\gamma} +\Theta_{\gamma}-\Theta_{\beta} }\\
&\leq \norm{ \Theta_{\alpha} -\Theta_{\gamma} } + \norm{ \Theta_{\gamma}-\Theta_{\beta}}.
\end{align*}
\goodbreak Let us denote the elements in the partitions $\alpha$ and $\gamma$ by 
\par
	  \noindent\mbox{
\begin{minipage}{0.50\textwidth}
\begin{align*}
\gamma =& \left\{ R=t_{1}<\dotsc < t_{n+1} =S  \right\} \\
\alpha =&\left\{  t_{1}=s_{1}^{(1)}<\dotsc < s_{m_{1}}^{(1)}=t_{2} \right\} \\
&\cup \left\{  t_{2}=s_{1}^{(2)}<\dotsc < s_{m_{2}}^{(2)}=t_{3} \right\} \\
&\dotsb \text{ }\\
 &\cup \left\{  t_{n}=s_{1}^{(n)}<\dotsc < s_{m_{n}}^{(n)}=t_{n+1}\! \right\}\!.
\end{align*}
\end{minipage}
}
	  \noindent\mbox{
\begin{minipage}{0.45\textwidth}
	     \includegraphics[height=4cm]{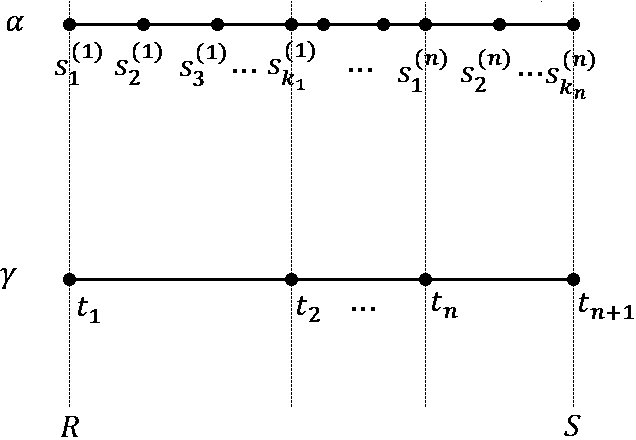}
\end{minipage}
}

\par
 Since $\alpha\geq \gamma$ we can define
\begin{align*}
\alpha_{i}:= \alpha \cap [t_{i},t_{i+1}] = \left\{  t_{i}=s_{1}^{(i)}<\dotsc < s_{m_{i}}^{(i)}=t_{i+1} \right\}
\end{align*}
for $i=1 \ldots n$. In particular, $\alpha_{i} \in \P{[t_{i},t_{i+1}]}$. Consider the factors  of $\Theta_{\alpha}$ between $t_{i}=s_{1}^{(i)}$ and $s_{m_{i}}^{(i)}=t_{i+1}$, i.e.
\begin{align*}
d_{i}:=& g_{s_{1}^{(i)},s_{2}^{(i)}}\dotsb g_{s_{m_{i}-1}^{(i)},s_{m_{i}}^{(i)}} %(1+a_{s_{1}^{i},s_{2}^{i}}+b_{s_{1}^{i},s_{2}^{i}} ) \dotsb (1+ a_{s_{m_{i}-1}^{i},s_{m_{i}+1}^{i}}+b_{s_{m_{i}-1}^{i},s_{m_{i}+1}^{i}})
\end{align*}
and the factor $g_{t_{i},t_{i+1}}$ of $\Theta_{\gamma}$, i.e.
\begin{align*}
c_{i}:=& g_{t_{i},t_{i+1}} \overset{\text{(\ref{PCondE})}}{=} 1+ \sum\limits_{k=1}^{n_{i}-1} (a_{s_{k}^{(i)},s_{k+1}^{(i)}})
\end{align*}
for $i=1, \ldots, n$. Then
\begin{align}
 &\norm{ \Theta_{\gamma} -\Theta_{\alpha} } \notag\\
 &=\norm{c_{1} \dotsb c_{n}  -d_{1}\dotsb d_{n} } \notag\\
 &=\norm{ \sum\limits_{j=1}^{n} c_{1} \dotsb c_{j-1} \Big( c_{j} - d_{j}\Big) d_{j+1}\dotsb d_{n} } \notag\\
 &\leq \sum\limits_{j=1}^{n} \norm{c_{1}} \dotsb \norm{c_{j-1}} \norm{ d_{j} - c_{j}}  \norm{d_{j+1}}  \dotsb \norm{d_{n} } \notag\\
 &\leq \sum\limits_{j=1}^{n} \expo{}{C(t_{j}-t_{1}+ t_{n+1}-t_{j+1})}{}  \norm{ d_{j} - c_{j}} \label{PGoOn}
\end{align}
with telescoping sum $c_{1} \dotsb c_{n}  -d_{1}\dotsb d_{n}=\sum\limits_{j=1}^{n} c_{1} \dotsb c_{j-1} \Big( c_{j} - d_{j}\Big) d_{j+1}\dotsb d_{n} $.
Applying Lemma \ref{BehII} to $\norm{d_{j} - c_{j}}$ yields
\begin{align*}
\norm{ d_{j} - c_{j}} \,\leq\,\, \left( \frac{1}{2}(t_{j+1}-t_{j})^2C^2\expo{}{C(t_{j+1}-t_{j})}{} \right),
\end{align*}
where  $(1+a_{1})\dotsb (1+a_{n}) := d_{j}$ and $1+a_{1} +\dotsb +a_{n} := c_{j}$.\par
Let $\norm{\gamma}  := \max \left\{  t_{i+1}-t_{i}, i= 1\dotsc n \right\}$ and continue at formular (\ref{PGoOn}), so
\begin{align*}
\text{(\ref{PGoOn})}&\leq \frac{1}{2}C^2 \expo{}{C(S-R)}{} \sum\limits_{j=1}^{n} (t_{j+1}-t_{j}) \norm{\gamma}\\
&= \frac{1}{2}C^2 \expo{}{C(S-R)}{}\norm{\gamma}(S-R).
\end{align*}
Therefore, 
\begin{align*}
\norm{ \Theta_{\alpha} -\Theta_{\gamma} } + \norm{ \Theta_{\gamma}-\Theta_{\beta} } \leq \norm{\gamma} C^2 \expo{}{C(S-R)}{}(S-R).
\end{align*}
This tends to $0$ if $\norm{\gamma} \to 0$.\par
It remains to show that for every $\epsilon>0$ there exists a partition $\gamma$ such that $ \norm{\gamma}C^2 \expo{}{C(S-R)}{}(S-R) < \epsilon$. Choosing a sufficiently fine equidistant partition $\gamma$ finishes the proof. 
\end{beweis}
\section{Applications}\label{sec:example}
In this section, we apply the previous main lemma to example I,II and III.
Example I is a special case of example II for $d=1$. Both, example II and III originate from the same construction namely convolving a familiy of algebra homomorphisms with respect to the convolution $\astalg$ given by the comultiplication $\Delta$ of a dual semigroup. In order to see that the matrices from example III  arise from this kind of construction, we consider a family of algebra homomorphisms $f_{r,s}$ together with the primitve comultiplication $\Delta_{p}$ of \ref{dSG:tensoralgebra} and denote the kernel of the counit $\delta$ by $B_{0}$:
\begin{align*}
 \left( \bigop{\star_{\Delta_{p}}}{i=1}{n} f_{t_{i},t_{i+1}} \right) (\delta_{k,l}\eins+ x_{k,l}-\delta_{k,l}\eins) &= \delta_{k,l}\eins+\sum\limits_{i=1}^{n} f_{t_{i},t_{i+1}}(x_{k,l}-\delta_{k,l}\eins) \\
 &= (1-n)\delta_{k,l}\eins+\sum\limits_{i=1}^{n} f_{t_{i},t_{i+1}}(x_{k,l}). 
\end{align*}
We invite the reader to compare the last expression of this equation to the entries of the matrices the net in example III consists of.
In example II, we consider the dual semigroup $\KD $ from \ref{dSG:kd}.
\subsection{Example I}
Let $A=B(Fock(L_{2}(\R_{+})  \otimes \C ))$. Define
\begin{align*}
a_{r,s}:=& A^{*}(\chi_{\left[ r,s \right]}  \otimes 1)- A(\chi_{\left[ r,s \right]}   \otimes 1) -\frac{1}{2}(r-s)id,
\end{align*}
where $A^{*}$ and $A$ denote the creation and the annihilation operator on the full Fock space and $\chi_{\left[ r,s \right]}$ the characteristic function of the intervall $\left[ r,s \right]$. Since the family
 $(a_{r,s})_{0\leq r<s}$ satisfies conditions (\ref{PCondC}) and (\ref{PCondE}) Lemma \ref{abstr:sz:cauchyNet} applies to the net 
\begin{align*}
(\Theta_{\alpha})_{\alpha \in \P{[R,S]}}:= \left( (id + a_{t_{1},t_{2}} )\dotsb ( id + a_{t_{n},t_{n+1}} ) ) \right)_{\alpha=\{t_{1}<\dotsb < t_{n+1}\} \in \P{[R,S]}}.
\end{align*}
Therefore, this net is a Cauchy net.

\subsection{Example II}\label{sec:ugroup}
Consider a dual semigroup $(B,\Delta,\delta)$. 
\begin{df} A generator $\Psi$ on $\pa{ B,\delta }$ is a linear map $\Psi:B\to \C$ with\goodbreak
$\Psi\left( b^{*} \right)=\overline{\Psi\left( b \right)}$ for all $b\in B$, $\Psi\left( \eins{} \right)=0$ and  $\Psi\left( b^{*}b \right)\geq 0$ for all $b\in ker\left( \delta \right)$.
\end{df}
\begin{df}[Sch\"urmann triple] \label{def:uc}
A Sch\"urmann triple $\pa{ \rho,\eta,\Psi }$ of $\pa{ B,\delta }$ on a pre-Hilbert space  $D$ consists of 
\begin{itemize}
\item[(i)] a unital $*$\nbd algebra homomorphism $\rho:B\to L_{a}(\overline{D})$ ($*$\nbd representation) 
\item[(ii)] a surjective, linear map $\eta:B\to D$ with
\begin{align}
\eta\pa{ ab }= \rho\pa{ a }\eta\pa{ b }+ \eta\pa{ a }\delta\pa{ b }
\end{align}
\item[(iii)] a generator $\Psi:B\to \C$ with
\begin{align}
\pd{ \eta\pa{ a },\eta\pa{ b } }_{D}=\Psi\pa{ \pa{ a-\delta\pa{ a }\eins }^{*}\pa{ b-\delta\pa{ b }\eins } }.
\end{align}
\end{itemize}
\end{df}
There exists a construction of a Sch\"urmann triple of $\Psi$ which resembles the GNS-construction \cite[Prop 4.1]{MSchue90a}. Let $B$ be the dual semigroup $\KD $  for $d\geq 1$ from example \ref{dSG:kd}. The existence of a generator $\Psi$ on $\KD $ is ensured by \cite[Theorem 5.1.12]{MSchue93}:
\begin{sz}\label{kdGenerator}Let $L$ be a complex $d\times d$-matrix for $d\geq1$. Then there exists a unique generator $\Psi:\KD  \to \C$ and a Sch\"urmann triple $\pa{ \rho,\eta,\Psi }$ such that 
$\Psi\left( x_{kl} \right)= \frac{1}{2}\left( L^{*}L \right)_{kl}$, $D = \C$, $\eta\left( x_{kl} \right) = L _{kl}$, $\eta\left( x_{kl}^{*} \right)= -  L _{lk}$ and $\rho\left( x_{kl} \right) = \delta_{kl}$ for all $k,l = 1,\ldots, d$.
\end{sz}

Let $\Psi:B\to \C$ be a generator with Sch\"urmann triple $\pa{ \rho,\eta,\Psi }$.  Consider the full Fock space $\Gamma$ over $L_{2}\left( \R_{+},D \right)\cong D   \otimes L_{2}\left( \R_{+} \right)$ and
% \begin{align*}
% L_{2}\left( \R_{+},D \right)\cong D   \otimes L_{2}\left( \R_{+} \right)\cong D \otimes L_{2}\left( \left[ S,T \right] \right)  \oplus D \otimes L_{2}\left( \left[ S,T \right]^{c} \right)
% \end{align*}
define a family of unital $*$\nbd algebra homomorphisms $h_{r,s}:B \to L_{a}(\Gamma)$ by
\begin{align}\label{LP:hrs}
h_{r,s}(x_{k,l}):=\delta_{k,l}id+ A_{r,s}^{*} (\eta (x_{k,l}))  +A_{r,s} (\eta ((x_{k,l})^{*})) +(s-r)\Psi(x_{k,l})id,
\end{align}
where $A^{*}_{r,s}(d):= A^{*}(\chi_{\left[ r,s \right]}  \otimes d)$ and $A_{r,s}(d):= A(\chi_{\left[ r,s \right]}  \otimes d)$.
Let $a_{r,s}$ be the matrix with entries 
\begin{align}
(a_{r,s})_{k,l}:=h_{r,s}(x_{k,l}-\delta_{k,l}\eins).
\end{align}
The family $(a_{r,s})_{0\leq r<s}$ satisfies conditions (\ref{PCondC}) and (\ref{PCondE}). Applying Lemma \ref{abstr:sz:cauchyNet} implies that the net 
\begin{align*}
(\Theta_{\alpha})_{\alpha \in \P{[R,S]}}:= \Bigl(  \prod\limits_{\alpha} g_\alpha\Bigr)_{\alpha \in \P{[R,S]}}%(id + a_{t_{1},t_{2}} )\dotsb ( id + a_{t_{n},t_{n+1}} ) ) \right)_{\alpha=\{t_{1}<\dotsb < t_{n+1}\} \in \P{[R,S]}}
\end{align*}
is a Cauchy net, where $g_{r,s}:=id+a_{r,s}$. Furthermore, the elements of the matrix $g_{r,s}$  are
\begin{align*}
(g_{r,s})_{k,l}= h_{r,s}((E+(X-E))_{k,l})=h_{r,s}(x_{k,l}).
\end{align*}
For a partition $\alpha=\{t_{1}<\dotsb < t_{n+1}\} \in \P{[R,S]}$ we get that
\begin{align}
\left(\Theta_{\alpha}\right)_{k,l}&= (g_{t_{1},t_{2}}\dotsm g_{t_{n},t_{n+1}})_{k,l} \nonumber\\
&=\left( h_{t_{1},t_{2}} \sqcup_{1} \dotsb  \sqcup_{1} h_{t_{n},t_{n+1}}  \right) \left( (X^{(1)}\dotsm X^{(n)})_{k,l} \right) \nonumber\\
&=\left( h_{t_{1},t_{2}} \sqcup_{1} \dotsb  \sqcup_{1} h_{t_{n},t_{n+1}}  \right)  \circ \Delta_{n}(x_{k,l}) \nonumber\\
&=\bigl(  \bigop{ \sqcup_{1}}{\alpha}{} h_{\alpha}\bigr)  \circ \Delta_{n}(x_{k,l}) \label{NetAlg},
\end{align}
where $\Delta_{2}:=\Delta$ and $\Delta_{n+1}:=\left( \Delta_{n} \sqcup_{1} id \right)  \circ \Delta$. Therefore, the matrix entry $\left(\Theta_{\alpha}\right)_{k,l}$ is the evaluation of a unital $*$\nbd algebra homomorphism at $x_{k,l}$.
\begin{df}
Let $U_{R,S}$ be the limit of the net $(\Theta_{\alpha})_{\alpha \in \P{[R,S]}}$.
\end{df}
The net $(\Theta_{\alpha})_{\alpha \in \P{[R,S]}}$ and its limit $U_{R,S}$ are $d\times d$ matrices with entries in the algebra $L_{a}(\Gamma (D   \otimes L_{2}\left( \R_{+} \right)))$.
\begin{sz} The family $\pa{ U_{R,S} }_{0\leq R\leq S}$ has the following properties: 
\begin{itemize}
\item $U_{r,s}$ is a unitary matrix for $0\leq r\leq s$
\item $U_{r,s}U_{s,t}=U_{r,t}$ for $0\leq r<s<t$
\item $\lim\limits_{r<s, s\to r} U_{r,s}= E$ for $0\leq r$.
\end{itemize}

\end{sz}
\begin{beweis} Since $U_{R,S}$ is the norm limit of $(\Theta_{\alpha})_{\alpha \in \P{[R,S]}}$, we have
\begin{align*}
(U_{R,S}^{*}U_{R,S})_{k,l} =& \lim\limits_{\alpha \in \P{[R,S]}, \norm{\alpha} \to 0} (\Theta_{\alpha}^{*}\Theta_{\alpha})_{k,l}.
\end{align*}
Consider a partition $\alpha=\{t_{1}<\dotsb < t_{n+1}\} \in \P{[R,S]}$. Then
\begin{align*}
&(\Theta_{\alpha}^{*}\Theta_{\alpha})_{k,l}= \sum\limits_{p=1}^{d} \left( \Theta_{\alpha}^{*} \right)_{k,p} \left( \Theta_{\alpha} \right)_{p,l} \\
&=\sum\limits_{p=1}^{d}  \left(  \bigop{\sqcup_{1}}{\alpha}{} h_{\alpha} \right)  \circ \Delta_{n}(x_{k,p}^{*})  \cdot 
 \left(  \bigop{\sqcup_{1}}{\alpha}{} h_{\alpha} \right)  \circ \Delta_{n}(x_{p,l}) \\
% &=\sum\limits_{p=1}^{d} \!\! \left( h_{t_{1},t_{2}} \sqcup_{1} \dotsb  \sqcup_{1} h_{t_{n},t_{n+1}}  \right)  \circ \Delta_{n}(x_{k,p}^{*})  \cdot 
% \left( h_{t_{1},t_{2}} \sqcup_{1} \dotsb  \sqcup_{1} h_{t_{n},t_{n+1}}  \right)  \circ \Delta_{n}(x_{p,l}) \\
&= \left(  \bigop{\sqcup_{1}}{\alpha}{} h_{\alpha} \right)  \circ \Delta_{n}(\underbrace{\sum\limits_{p=1}^{d}x_{k,p}^{*}x_{p,l}}_{=(X^{*} \cdot X)_{k,l}=\delta_{k,l}\eins}) = \delta_{k,l}id 
\end{align*}
which implies $U_{R,S}^{*}U_{R,S}= E$. Analogously, $U_{R,S}U_{R,S}^{*} =E $. For $1\leq k,l \leq d$ it holds:
\begin{align*}
 \left( U_{r,s} \cdot U_{s,t} \right)_{k,l} &= \Bigl( \lim\limits_{\substack{\alpha=\alpha_{1}\cup\alpha_{2}\\
\alpha_{1}=\left\{  r=s_{1}<\dotsb < s_{m+1}=s \right\}\\
\alpha_{2}=\left\{  s=t_{1}<\dotsb < t_{n+1}=t \right\}}} \underbrace{\Theta_{r,s,\alpha_{1}}\Theta_{s,t,\alpha_{2}}}_{=\Theta_{r,t,\alpha}} \Bigr)_{k,l} =\left( U_{r,t}  \right)_{k,l}.
\end{align*}
When withdrawing the term $- \sum\limits_{i=1}^{n} a_{i}$ in Lemma \ref{BehII} modifying the proof yields the norm limit $\lim\limits_{r<s, s\to r} U_{r,s}= E$.
\end{beweis}
The $h_{r,s}$ from (\ref{LP:hrs}) consist of freely adapted creation and annihilation operators. Thus, $g_{r,s}$ are freely adapted entry-wise. Therefore, the limit $U_{R,S}$ is freely adapted at $D   \otimes L_{2}\left( \pb{ R,S } \right)$ entry-wise. In order to present a unitary cocycle equation we can replace $\R_{+}$ by $\R$ in our consideration, especially for  $L_{2}\left( \R \right)$ instead of  $L_{2}\left( \R_{+} \right)$ and use time shifts $\mathfrak{s}_{r}$ in $L_{a}(\Gamma(D   \otimes L_{2}\left( \R \right)))$, where 
\begin{align*}
\mathfrak{s}_{r}(x_{1} \otimes \dotsb  \otimes x_{p})= (\mathfrak{s}_{r}'(x_{1})) \otimes \dotsb  \otimes (\mathfrak{s}_{r}'(x_{p}))
\end{align*}
and $\mathfrak{s}_{r}'(x) (t)=x(t+r)$ denote the time shifts in $L_{2}\left( \R ,D\right)$. When defining\goodbreak  $W_{t}:=U_{0,t}$ the following unitary cocycle equation holds:
\begin{align*}
W_{s+t}= W_{t}\,\mathfrak{s}_{t}\,W_{s}\,\mathfrak{s}^{*}_{t}
\end{align*}
\subsection{Example III}\label{sec:tensoralg}
Conversely, the family of unitaries $U_{r,s}$ yields a net converging to $g_{R,S}$.
\begin{df}
For $\alpha=\{R=t_{1}<\dotsb < t_{n+1}=S\} \in \P{[R,S]}$ define
\begin{align*}
\Upsilon_{\alpha}:=(1-n)id+ \sum\limits_{j=1}^{n} U_{t_{j},t_{j+1}}.
\end{align*}

\end{df}

\begin{sz}
The net $(\Upsilon_{\alpha})_{\alpha \in \P{[R,S]}}$ 
%$\left((1-n)id+ \sum\limits_{j=1}^{n} U_{t_{j},t_{j+1}} \right)_{\alpha=\{R=t_{1}<\dotsb < t_{n+1}=S\} \in \P{[R,S]}}$\par  
converges to $g_{R,S}$.
\end{sz}
\begin{beweis}
Let $\alpha=\{R=t_{1}<\dotsb < t_{n+1}=S\} \in \P{[R,S]}$ and $\epsilon>0$. By (\ref{PCondE}) 
\begin{align*}
(1-n)id-g_{R,S} =- \sum\limits_{j=1}^{n} g_{t_{j},t_{j+1}} 
\end{align*}
holds.
For $i=1, \ldots, n$ a partition $\beta_{i}\in\P{[t_{i},t_{i+1}]}$ with 
\begin{align*}
\norm{  \prod_{\beta_{i}} g_{\beta_{i}} -U_{t_{i},t_{i+1}} } < \frac{\epsilon}{n}
\end{align*}
exists due to the convergence of $(\Theta_{\alpha_{i}})_{\alpha_{i}\in\P{[t_{i},t_{i+1}]}}$ to $U_{t_{i},t_{i+1}}$ in example II. 
Lemma \ref{BehII} yields
\begin{align*}
\norm{\prod_{\beta_{i}} g_{\beta_{i}} - g_{t_{i},t_{i+1}}} \overset{\text{(\ref{BehII})}}{\leq}  \frac{1}{2}(t_{i+1}-t_{i})^2C^2 \expo{}{C(t_{i+1}-t_{i})}{}.
\end{align*}
% For $u=t_{j},v=t_{j+1}$ a partition $\beta_{u,v}=\{u=s_{1}<\dotsb < s_{m+1}=v\}$ with 
% \begin{align*}
% \norm{ g_{s_{1},s_{2}}\dotsb g_{s_{m},s_{m+1}} -U_{u,v} } < \frac{\epsilon}{n}
% \end{align*}
% exists for all $j=1,\ldots, n$ due to the convergence of $(\Theta_{\alpha})_{\alpha\in\P{[u,v]}}$ to $U_{u,v}$ in example II. 
% Lemma \ref{BehII} yields
% \begin{align*}
% \norm{g_{s_{1},s_{2}}\dotsb g_{s_{m},s_{m+1}} - g_{u,v}} =& \norm{(1+a_{1})\dotsb (1+a_{m}) -1 -a_{1}-\dotsb - a_{m}}\\
% \overset{\text{(\ref{BehII})}}{\leq}&  \frac{1}{2}(v-u)^2C^2 \expo{}{C(v-u)}{}
% \end{align*}
% with $a_{j}=a_{s_{j},s_{j+1}}$ and $g_{u,v}\overset{\text{(\ref{PCondE})}}{=} 1+a_{1}+\dotsb + a_{m}$.
% Therefore,
% \begin{align*}
% &\norm{(1-n)id+\sum\limits_{j=1}^{n} U_{t_{j},t_{j+1}}-g_{R,S}} \\
% =& \norm{\sum\limits_{j=1}^{n} \left( U_{t_{j},t_{j+1}} \pm ( g_{s_{1},s_{2}}\dotsb g_{s_{m},s_{m+1}})_{s_{1}:=t_{j}, s_{m+1}:=t_{j+1}} -g_{t_{i},t_{i+1}}\right)  } \\
% \leq& \epsilon + \sum\limits_{j=1}^{n}\norm{(g_{s_{1},s_{2}} \dotsb g_{s_{m},s_{m+1}})_{s_{1}:=t_{j}, s_{m+1}:=t_{j+1}} - g_{t_{j},t_{j+1}}} \\
% \leq& \epsilon + \sum\limits_{j=1}^{n} \frac{1}{2}(t_{j+1}-t_{j})^2C^2 \expo{}{C(t_{j+1}-t_{j})}{}\\
% \leq& \epsilon +\frac{1}{2} \norm{\alpha}(S-R)C^2\expo{}{C(S-R)}{}.
% \end{align*}
Therefore,
\begin{align*}
&\norm{(1-n)id+\sum\limits_{j=1}^{n} U_{t_{j},t_{j+1}}-g_{R,S}} \\
&= \norm{\sum\limits_{j=1}^{n} \Bigl( U_{t_{j},t_{j+1}} \pm \Bigl(\prod_{\beta_{j}} g_{\beta_{j}} \Bigr) -g_{t_{j},t_{j+1}}\Bigr)  } \\
&\leq \epsilon + \sum\limits_{j=1}^{n}\norm{\prod_{\beta_{j}} g_{\beta_{j}} - g_{t_{j},t_{j+1}}} \\
&\leq \epsilon + \sum\limits_{j=1}^{n} \frac{1}{2}(t_{j+1}-t_{j})^2C^2 \expo{}{C(t_{j+1}-t_{j})}{}\\
&\leq \epsilon +\frac{1}{2} \norm{\alpha}(S-R)C^2\expo{}{C(S-R)}{}.
\end{align*}
It remains to show that for $\delta>0$ there exists a partition $\alpha$ with\goodbreak
$ \frac{1}{2} \norm{\alpha}(S-R)C^2\expo{}{C(S-R)}{} <\frac{\delta}{2}$. Let $\epsilon<\frac{\delta}{2}$ and choose a sufficiently fine equidistant partition $\alpha$. This finishes the proof. 
\end{beweis}

\section{Application to free L\'evy processes}\label{sec:QLP}
This section connects the previous sections by using the families $(h_{r,s})_{0\leq r\leq s}$ and $\big(U_{R,S}\big)_{{0\leq R\leq S}}$ from example II and III to construct a free L\'evy process over $\KD $ on the full Fock space.
We define  free (quantum) L\'evy processes over a dual semigroup on a quantum probability space.
\begin{df}[L\'evy process]
Let $(B,\Delta,\delta)$ be a dual semigroup. A family $\left( f_{r,s} \right)_{0\leq r\leq s}$ of unital $*$\nbd algebra homomorphisms $f_{r,s}:(B,\Delta,\delta)\to (A,\Phi)$ with a state $\Phi:A\to \C$ such that
\begin{itemize}
\item[  (i)] $r<s<t \Rightarrow f_{r,t}=f_{r,s}\astalg f_{s,t}$ and $f_{r,r}(b)=\delta(b)id \,\,\forall b \in B$
\item[ (ii)]  $n>1\in \N, 0\leq t_{1}<\dotsb <t_{n+1} \Rightarrow$ 
\begin{align*}
\Phi \circ \left( f_{t_{1},t_{2}} \sqcup \dotsb  \sqcup f_{t_{n},t_{n+1} }\right)( \cdot)=\left( \Phi \circ f_{t_{1},t_{2}} \right) \odot \dotsb  \odot \left( \Phi \circ f_{t_{n},t_{n+1}}\right) ( \cdot)
\end{align*}
\item[(iii)] $r<s \Rightarrow \Phi \circ f_{r,s}=\Phi \circ f_{0,s-r}$
\item[(iv)] $b\in B \Rightarrow \lim\limits_{t\to 0} \Phi \circ f_{0,t}(b)=\delta(b)$
\end{itemize}
is called a free Lévy process.  
\end{df}
\begin{bm0}
For an arbitary $\odot$ from Definition \ref{PUP} and (UP4) the definition is the same.
\end{bm0}

\begin{sz}\label{MarginalDistribution}
The family of marginal distributions $\varphi_{t}:= \Phi \circ f_{0,t}$ for $t\geq 0$ is a (weakly) continuous convolution (w.r.t $\uplus$) semigroup of states. The derivation $\lim\limits_{t\to 0} \frac{1}{t}(\varphi_{t}(b)-\delta(b)) = \Psi(b)$ exists and defines a generator $\Psi$. 
\end{sz}
\begin{beweis}
See \cite[Prop4.4, Theorem4.6]{BGhSc}.
\end{beweis}
\begin{sz}[additive free white noise]\label{addWhiteNoise}
Consider the dual semigroup $T\left( V \right)$ from example \ref{dSG:tensoralgebra}. For a generator $\Psi:T(V) \to \C$ and a Sch\"urmann triple $\pa{ \rho,\eta,\Psi }$ we define a unital $*$\nbd algebra homomorphism $T(I_{s,t}):T(V) \to L_{a}(\Gamma)$ by the  linear map
\begin{align*}
I_{s,t}(v):= A_{s,t}^{*} \circ \eta (v) + \Lambda_{s,t} \left( \rho(v) \right) + A_{s,t} \circ \eta (v^{*}) + (t-s)\Psi(v)id \quad \forall v\in V.
\end{align*}
The family $\left( T(I_{r,s}) \right)_{0\leq r \leq s}$ is a cyclic Lévy process with respect to freeness and $\Psi$ is the generator of the marginal distribution. Moreover, it is cyclic for $\C\eins  \oplus V$.
\end{sz}
\begin{beweis}
See \cite{GlScSp} and \cite{Fra06}.
\end{beweis}

\begin{bsp}
We now apply Theorem \ref{addWhiteNoise} to $\left(\KD , \Delta,\delta\right)$:\par
 Let $C:=ker(\delta)$. Then $B:=\KD  \cong \C\eins   \oplus C$ and $B \sqcup_{1} B \cong \C\eins  \oplus C \sqcup C$. Since $\delta(x_{k,l})=\delta_{k,l}$ the elements $x_{k,l}^{z}-\delta_{k,l}\eins$ are in the kernel of $\delta$. We now eva\-lu\-ate $T(I_{s,t}):T(C) \to L_{a}(\Gamma)$ at $x_{k,l} -\delta_{k,l}\eins$: 
\begin{align*}
I_{s,t}(x_{k,l} -\delta_{k,l}\eins):=& A_{s,t}^{*} \circ \eta (x_{k,l} -\delta_{k,l}\eins) + \Lambda_{s,t} \left( \rho(x_{k,l} -\delta_{k,l}\eins) \right) \\
&+ A_{s,t} \circ \eta ((x_{k,l} -\delta_{k,l}\eins)^{*}) + (t-s)\Psi(x_{k,l} -\delta_{k,l}\eins)id\\
=&A_{s,t}^{*} \circ \eta (x_{k,l} ) + \Lambda_{s,t} \left( \rho(x_{k,l}) -\delta_{k,l}\eins \right) \\
&+ A_{s,t} \circ \eta ((x_{k,l})^{*}) + (t-s)\Psi(x_{k,l})id,
\end{align*}
since $\eta(\eins)=0=\Psi(\eins)$.  If $\rho(x_{k,l}) =\delta_{k,l}\eins$ then
\begin{align*}
&T(I_{r,s})\big(\delta_{k,l}\eins +  (x_{k,l}-\delta_{k,l}\eins) \big) = h_{r,s}(x_{k,l}). 
\end{align*}
We conclude this example by stressing that the $h_{r,s}$ from example II are free.
\end{bsp}
For a generator and a Sch\"urmann triple induced by a complex $d\times d$ matrix $L$ as in Theorem \ref{kdGenerator} we know that $\rho(x_{k,l}) =\delta_{k,l}\eins$. 
\subsubsection*{Construction of a Lévy process over $\KD $}
Let $\Psi:\KD :\to \C$ be generator and let $\pa{ \rho,\eta,\Psi }$ be a Sch\"urmann triple with 
\begin{align}
\rho(x_{k,l}) =\delta_{k,l}\eins  \quad \text{ for } k,l=1, \ldots,d. \label{Gen} 
\end{align}
Define a family of unital $*$\nbd algebra homomorphisms by
\begin{align*}
&f_{r,s}:\C\langle d \rangle  \to (L_{a}(\Gamma),\langle\Omega, (\cdot) \Omega\rangle)\\
&f_{r,s}\left(x_{k,l}\right) =\left(U_{r,s}\right)_{k,l}
\end{align*}
and $f_{r,r}(b):=\delta(b)id$. Due to the unitaries $U_{r,s}$ the $f_{r,s}$ respect the relations of $\KD $ and we get a family of unital $*$\nbd algebra homomorphisms $(f_{r,s})_{r\leq s}$ on $\KD $.
\begin{lm}
For all $b\in \KD $
\begin{align}
f_{R,S}(b)=\lim\limits_{\substack{\alpha \in \P{[R,S]} \\ \alpha=\{R=t_{1}<\dotsb < t_{n+1}=S\}}} \left( h_{t_{1},t_{2}} \sqcup_{1} \dotsb  \sqcup_{1} h_{t_{n},t_{n+1}}  \right)  \circ \Delta_{n} (b) \label{QLP:NetAlg}
\end{align}
holds.
\end{lm}
\begin{beweis}
Since there are only unital $*$\nbd algebra homomorphisms involved, it suffices to prove the claim for monomials $M:= x_{1}^{z_{1}}\dotsm x_{m}^{z_{m}}$, where $z_{j} \in \N$ and $x_{j}=x_{k_{j},l_{j}}$ or $x_{j}=x_{k_{j},l_{j}}^{*}$ for some $0<k_{j},l_{j}\leq d$. We  then use (\ref{NetAlg}) to obtain the following equation: 
\begin{align*}
f_{R,S}(M)=(f_{R,S}(x_{1}))^{z_{1}}\dotsm (f_{R,S}(x_{k}))^{z_{k}} = \prod\limits_{j=1}^{k} (f_{R,S}(x_{j}))^{z_{j}} \\
=\lim\limits_{\substack{\alpha \in \P{[R,S]} \\ \alpha=\{R=t_{1}<\dotsb < t_{n+1}=S\}}}\prod\limits_{j=1}^{k}  \left( \left( h_{t_{1},t_{2}} \sqcup_{1} \dotsb  \sqcup_{1} h_{t_{n},t_{n+1}}  \right)  \circ \Delta_{n} (x_{j})  \right)^{z_{j}} \\
=\lim\limits_{\substack{\alpha \in \P{[R,S]} \\ \alpha=\{R=t_{1}<\dotsb < t_{n+1}=S\}}}  \bigg( \left( h_{t_{1},t_{2}} \sqcup_{1} \dotsb  \sqcup_{1} h_{t_{n},t_{n+1}}  \right)  \circ \Delta_{n} (\underbrace{\prod\limits_{j=1}^{k} x_{j}^{z_{j}}}_{=M})  \bigg) .
\end{align*}
\end{beweis}
\begin{sz}
The family $\left( f_{r,s} \right)_{0\leq r \leq s}$ is a Lévy process with respect to freeness and $\Psi$ is the generator of the marginal distribution.
\end{sz}
\begin{beweis}
\begin{tabbing}
(i) Evolution: \= $
f_{r,s}\astalg f_{s,t}  (x_{k,l})= f_{r,s} \sqcup f_{s,t} \left( X^{(1)} X^{(2)} \right)_{k,l}$\\
{} \> $=f_{r,s} \sqcup f_{s,t} \biggl(  \sum\limits_{n=1}^{d} x_{k,n}^{(1)} x_{n,l}^{(2)} \biggr) = \sum\limits_{n=1}^{d} f_{r,s}(x_{k,n}) f_{s,t}(x_{n,l})$ \\
{} \> $=\left( U_{r,s} \cdot U_{s,t} \right)_{k,l} = \left( U_{r,t}  \right)_{k,l} = f_{r,t}(x_{k,l})$.
\end{tabbing}
%
%
% (i) Evolution:
% \begin{align*}
% &f_{r,s}\astalg f_{s,t}  (x_{k,l})= f_{r,s} \sqcup f_{s,t} \left( X^{(1)} X^{(2)} \right)_{k,l} \\
% &=f_{r,s} \sqcup f_{s,t} \biggl(  \sum\limits_{n=1}^{d} x_{k,n}^{(1)} x_{n,l}^{(2)} \biggr) = \sum\limits_{n=1}^{d} f_{r,s}(x_{k,n}) f_{s,t}(x_{n,l}) \\
% &= \left( U_{r,s} \cdot U_{s,t} \right)_{k,l} = \left( U_{r,t}  \right)_{k,l} = f_{r,t}(x_{k,l}).
% \end{align*}
(ii) Freeness: The continuity of the scalar product $\langle \cdot, \cdot\rangle$ implies that the net\goodbreak $\Bigl( \Phi \circ  \bigop{ \sqcup_{1}}{i=1}{n}h_{\alpha_{i}} \Bigr) (a_{1} \dotsm a_{m})$ converges to $\Bigl( \Phi \circ  \bigop{ \sqcup_{1}}{i=1}{n}f_{t_{i},t_{i+1}} \Bigr) (a_{1} \dotsm a_{m})$ with partition\goodbreak $\alpha=\alpha_{1}\cup \ldots\cup \alpha_{n}\in \P{I_{1}\cup \ldots \cup I_{n}}$. The free $\odot$-product evaluated at a point is a polynomial, see (\ref{fp::A}). Due to the continuity of polynomials the net \goodbreak $\Bigl(   \bigop{ \ \odot}{i=1}{n}\Phi \circ h_{\alpha_{i}} \Bigr) (a_{1} \dotsm a_{m})$  converges to $\Bigl(   \bigop{  \odot}{i=1}{n}\Phi \circ f_{t_{i},t_{i+1}} \Bigr) (a_{1} \dotsm a_{m})$. Since the $h_{r,s}$ are free we get 
$%\begin{align*}
\Bigl( \Phi \circ  \bigop{ \sqcup_{1}}{i=1}{n}h_{\alpha_{i}} \Bigr) (a_{1} \dotsm a_{m})=\Bigl(   \bigop{ \ \odot}{i=1}{n}\Phi \circ h_{\alpha_{i}} \Bigr) (a_{1} \dotsm a_{m})
$.\par
(iii) Stationary: Use (\ref{QLP:NetAlg}), the freeness and the stationarity of $h_{r,s}$.\par
(iv) Weak continuity: Use (\ref{QLP:NetAlg}), the freeness, the stationarity and the weak continuity of $h_{r,s}$, then (UP3), the counit property and (UP1).\par
Generator: Observe that  
\begin{align*}
\Phi  \circ h_{t_{i},t_{i+1}} = \delta + (t_{i+1}-t_{i})\Psi+ R_{t_{i+1}-t_{i}} =: k_{t_{i+1}-t_{i}},
\end{align*}
 where $R_{t}$ and $k_{t}$ as in Theorem \ref{dualGroup:expo}. Then use (\ref{QLP:NetAlg}) and the freeness of the $h$'s to get the marginal distribution
\begin{align*}
\Phi \circ f_{0,t} &=\lim\limits_{\text{net}, \alpha} \Phi  \circ \bigop{ \star}{\alpha}{} h_{\alpha}
=   \lim\limits_{\text{net}, \alpha} \bigop{ \uplus}{\alpha}{} \Phi  \circ h_{\alpha}\\
&=  \lim\limits_{\text{net}, \alpha} \bigop{ \uplus}{\alpha}{} k_{\alpha}
=\exp_{\ast}(tD(\Psi)) 
\end{align*}
and 
$
\lim\limits_{t\to 0} \frac{1}{t}\bigl(-\delta+  \exp_{\ast}(tD(\Psi)) \bigr)(b)
=\lim\limits_{t\to 0} \frac{1}{t}\bigl(-\delta+ \delta +t\Psi +O(t^{2})\bigr)(b)
=\Psi(b)
$ yields the generator $\Psi$.
\end{beweis}
%
%
%
%
%\begin{bm0}
%Actually, assuming (UP1)-(UP3) from Definition \ref{PUP} it suffices to see that the $ \odot$-product evaluated at a point is a polynomial \cite{BGS02}.
%\end{bm0}
%
%

A L\'evy process on the full Fock space is called cyclic if  the evaluation of the unital $*$\nbd subalgebra generated by $f_{r,s}(b)$ for all $0\leq r\leq s$ and $b \in B$ at the vacuum vector $\Omega$ is dense in the full Fock space. From example III, we know that
$(1-n)\delta_{k,l}\eins+\sum\limits_{i=1}^{n} f_{t_{i},t_{i+1}}(x_{k,l})$ converges to $h_{R,S}(x_{k,l})$.  Since\goodbreak
$(1-n)\delta_{k,l}\eins =f_{0,1}((1-n)\delta_{k,l}\eins) $, $h_{R,S}(x_{k,l}) = T(I_{R,S})(x_{k,l})$ and the $T(I_{R,S})$ are cyclic over $\C\eins  \oplus B_{0}$, we get that 
\begin{fo}\label{co:cyclic}
The free L\'evy process $f_{r,s}$ is cyclic.
\end{fo}
By convolving $h_{r,s}$ (\ref{LP:hrs}) in example II, we have transferred the L\'evy process properties of the additive free white noise $T(I_{s,t})$ over a tensor algebra to unitaries $U_{R,S}$ which yields the free L\'evy process $f_{r,s}$ over $\KD $ with 'matrix multiplication' as comultiplication. Conversely, by convolving $f_{r,s}$ with respect to the additive comultiplication as in example III and by using the cyclic property of $T(I_{r,s})$  the L\'evy process $f_{r,s}$ is cyclic.

Two stochastic processes are called stochastically equivalent if their distributions coincide. By the stationarity, the marginal distribution determines the distribution of a L\'evy process and by Theorem \ref{MarginalDistribution} the marginal distribution is determined by a generator. 
We have shown a representation theorem for generators $\Psi:\KD \to \C$ with $\rho(x_{k,l}) =\delta_{k,l}\eins$ (\ref{Gen}), that is: 
A $\odot$-free L\'evy process over $\KD $ on any QPS with that generator is stochastically equivalent to the $\odot$-free L\'evy process $(f_{r,s})_{r\leq s}$ on the full Fock space from the previous construction.
\subsection*{Outline}
The methods established above expand beyond the setting of this paper. A general framework of computing  L\'evy processes over dual semigroups $(B,\Delta,\delta)$ using additive  L\'evy processes $T(I_{r,s})$ which also covers the case of $\rho(x_{k,l}) \not=\delta_{k,l}\eins$ and yields the representation theorem on Fock spaces not only in the case of freeness, is in preparation by the author.

\bibliographystyle{alpha}
\bibliography{mybib-1}

\end{document}